\documentclass[12pt]{article}
\usepackage{amsmath, amsthm, amsfonts, amssymb}
\usepackage{graphics}
\usepackage[font=small,format=plain,labelfont=bf,up]{caption}

\usepackage{hyperref}
\usepackage[english]{babel}

\usepackage{ifpdf} 
\ifpdf
\usepackage{epstopdf}
\fi

\theoremstyle{plain}

\newtheorem{theorem}{Theorem}

\theoremstyle{definition}

\newtheorem{example}[theorem]{Example}

\newcommand\ipic[1]		{\raisebox{-0.5\height}{\scalebox{.4}{\includegraphics{pic#1.eps}}}}
\newcommand\jpic[1]		{\raisebox{-0.5\height}{\scalebox{.5}{\includegraphics{pic#1.eps}}}}

\DeclareMathOperator{\Hom}{Hom}
\DeclareMathOperator{\Rep}{Rep}

\DeclareMathOperator{\Ad}{Ad}

\newcommand\vect{vect}
\newcommand\svect{svect}
\newcommand\id            {id}
\newcommand\eps           {\varepsilon}
\newcommand\one           {{\bf1}}
\newcommand\h	{\mathfrak{h}}

\newcommand\Cb            {\mathbb{C}}
\newcommand\Zb            {\mathbb{Z}}

\newcommand\Cc            {\mathcal{C}}

\newcommand\Sc            {\mathcal{S}}

\begin{document}

\def\thefootnote{\fnsymbol{footnote}}
\begin{flushright}
ZMP-HH/13-1\\
Hamburger Beitr\"age zur Mathematik 465
\end{flushright}

\vskip 1.5em

\begin{center}\Large
A braided monoidal category for symplectic fermions
\end{center}

\vskip 1.5em

\begin{center}
A.\ Davydov${}^1$, I.\ Runkel${}^2$\footnote{IR thanks the organisers of the {\em XXIX Internat.\ Colloq.\ on Group-Theoretical Methods in Physics} (August 20--26, 2012, Chern Inst., Tianjin) for the opportunity to speak.}
\\[.5em]
\emph{
${}^1$ Dept.\ Math., Ohio Univ., Athens, Ohio 45701, USA\\
${}^2$ Dept.\ Math., Hamburg Univ.,
Bundesstr.\ 55, 20146 Hamburg, Germany
}
\end{center}

\vskip 1.5em

\begin{abstract}
We describe a class of examples of braided monoidal categories which are built from Hopf algebras in symmetric categories. The construction is motivated by a calculation in two-dimensional conformal field theory and is tailored to contain the braided monoidal categories occurring in the study of the Ising model, their generalisation to Tamabara-Yamagami categories, and categories occurring for symplectic fermions.
\end{abstract}

\vskip 2em

\section{Introduction}\label{intro}

In this short note we summarise some of the results in \cite{Davydov:2012xg,Runkel:2012cf}, where also more extensive references can be found. 

We are interested in a particular type of $\Zb/2\Zb$-graded braided monoidal categories. The grade 0 component is the monoidal category $\Rep_{\Sc}(H)$ of modules over a Hopf algebra $H$ in a symmetric monoidal category $\Sc$, and the grade 1 component is the category $\Sc$ itself. We will write $\Cc = \Cc_0 + \Cc_1$ with $\Cc_0 = \Rep_{\Sc}(H)$ and $\Cc_1 = \Sc$. The tensor product functor $\ast$ on the various components is defined as:

\renewcommand{\arraystretch}{1.1}

\noindent
  \begin{tabular}{c@{\hspace{.6em}}c@{\hspace{1.0em}}l@{\hspace{5pt}}ll}
  $A$ & $B$ & $A \ast B$ && comments \\
  \hline
  $\Cc_0$ & $\Cc_0$ & $A \otimes B$ & $\in \Cc_0$ & the  $H$-action is via the coproduct of $H$
\\
  $\Cc_0$ & $\Cc_1$ & $F(A) \otimes B$ & $\in \Cc_1$ & $F : \Rep_{\Sc}(H) \to \Sc$ is the forgetful functor
\\
  $\Cc_1$ & $\Cc_0$ & $A \otimes F(B)$ & $\in \Cc_1$ & 
\\
  $\Cc_1$ & $\Cc_1$ & $H \otimes A \otimes B$ & $\in \Cc_0$ & the  $H$-action is by multiplication
  \end{tabular}

\bigskip
This somewhat ad-hoc looking definition of the tensor product is actually quite natural. The mixed tensor products are the natural left and right action of the monoidal category $\Rep_{\Sc}(H)$ on $\Sc$. To obtain the last line in the table, assume that $\Cc$ can be made rigid. Writing $T$ for the tensor unit of $\Sc$ considered as an object in $\Cc_1$, we have for all $H$-modules $M$
\begin{equation} \label{ty-tens}
  \Hom_{\Cc_0}(T^* \ast T , M) ~\cong~ \Hom_{\Cc_1}(T,T \ast M) ~\cong~ \Hom_{\Sc}(1,F(M)) \ .
\end{equation}
This means that $T^* \ast T$ is a representing object for the functor $M \mapsto \Hom_{\Sc}(1,F(M))$, and so $T^* \ast T \cong H$ as an $H$-module. If we in addition demand that $T^* \cong T$, the last line in the above table follows.

Given the above form of the tensor product functor $\ast$ on $\Cc$, one can ask if it is possible to describe associativity and braiding isomorphisms for $\ast$ in terms of Hopf algebraic data on $H$. Our results for this question are given in Section \ref{unif}. But before getting there, in Sections \ref{ty} and \ref{sf} we would like to give the two examples of such $\Zb/2\Zb$-graded braided monoidal categories which were our main motivation when setting up the formalism.

\section{Tambara-Yamagami categories}\label{ty}

For simplicity, we will work over the field $\Cb$.
Consider a fusion category $\Cc$ over $\Cb$ whose simple objects are labelled by $G \cup \{ m \}$ where $G$ is a finite group and $m$ is an extra label. Suppose that the tensor product $\ast$ is of the form, for $a,b \in G$,
\begin{equation} \label{ty-fusion}
  a \ast b \cong ab 
  ~~, \quad
  m \ast a \cong m \cong a \ast m
  ~~, \quad
  m \ast m \cong \textstyle \bigoplus_{g \in G} g \ .
\end{equation}
This tensor product is a special case of the one in the above table: the underlying symmetric category $\Sc$ is $\vect(\Cb)$, the category of finite dimensional $\Cb$ vector spaces. The  component $\Cc_0$ is spanned by the simple objects $g \in G$; the component $\Cc_1$ is spanned by $m$ alone, so that $\Cc_1 \cong \vect(k)$. The Hopf algebra $H \in \Sc$ is the function algebra $\mathrm{Fun}(G,\Cb)$. 

For any such fusion category $\Cc$, the group $G$ is necessarily abelian and $\Cc$ is monoidally equivalent to $\Cc(\chi,\tau)$, which is defined as follows \cite[Thm.\,3.2]{Tambara:1998}. $\Cc(\chi,\tau)$ has simple objects and fusion rules as in \eqref{ty-fusion}, and its associator is determined by a symmetric non-degenerate bicharacter $\chi : G \times G \to \Cb^\times$ and a choice of $\tau \in \Cb^\times$ such that $\tau^2 = |G|^{-1}$. The associator is a bit lengthy and we refer to \cite{Tambara:1998}.

The category $\Cc(\chi,\tau)$ allows for a braiding if and only if $G$ is an elementary 2-group (i.e.\ $gg=e$ for all $g \in G$) \cite{Siehler:2001}. The braiding isomorphisms are determined by a quadratic form $\sigma$ associated to the bicharacter\footnote{
This means that $\sigma : G \to \Cb^\times$ satisfies $\sigma(a) = \sigma(a^{-1})$, $\sigma(e)=1$, and that $\chi(a,b) \sigma(a) \sigma(b) = \sigma(ab)$ for all $a,b \in G$.} 
$\chi$ and a number $\beta$ such that $\beta\neq 0$ and $\beta^2 = \tau\,\sum_{a \in G} \sigma(a)$.
Explicitly, under the identifications \eqref{ty-tens} the braiding is ($a,b \in G$)
\begin{equation}\label{ty-braid}
 c_{a,b} = \chi(a,b) \, \id_{ab} 
 ~,~
 c_{a,m} = \sigma(a) \, \id_m = c_{m,a}
 ~,~
 c_{m,m} = \textstyle \bigoplus_{g \in G} \beta \,\sigma(g)^{-1} \, \id_g \ .
\end{equation}

An important example of a braided monoidal category of the above type is provided by the two-dimensional critical Ising model. There, one considers the three irreducible representations $\hat\one$, $\hat\eps$, $\hat\sigma$ of the Virasoro algebra which have central charge $c=\frac12$ and lowest $L_0$-weights $h_{\hat\one} = 0$, $h_{\hat\eps} = \frac12$ and $h_{\hat\sigma} = \frac{1}{16}$. The fusion rules are of the form \eqref{ty-tens} where $\hat\one,\hat\eps$ generate the group $G = \Zb/2\Zb$ and $m=\hat\sigma$ has fusion rule $\hat\sigma \ast \hat\sigma \cong \hat\one \oplus \hat\eps$. The braiding isomorphism $c_{r,s}$ -- projected to the simple object $t \in r \ast s$ -- is multiplication by $\exp(\pi i (h_r+h_s-h_t))$. Comparing to \eqref{ty-braid} shows that the braided monoidal structure is determined by
$\sigma(\hat\eps)= \exp(\pi i / 2)$, $\beta = \exp(\pi i / 8)$ and thus 
$\chi(\hat\eps,\hat\eps)=-1$, $\tau = 1/\sqrt{2}$.

\section{Symplectic fermions}\label{sf}

Continuing with examples from two-dimensional conformal field theory, we now consider symplectic fermions \cite{Kausch:1995py}. The mode algebra of $n$ pairs of symplectic fermsions is determined by a $2n$-dimensional symplectic vector space $\h$. It is convenient to think of $\h$ as a purely odd abelian Lie super-algebra with non-degenerate super-symmetric pairing $(-,-)$; we will use this language in the following. The symplectic fermion mode algebra is the affinisation $\hat\h$ of $\h$ with central element $K$ and graded bracket $[a_m,b_n] = m \,(a,b) \, \delta_{m+n,0} \, K$, where $m,n \in \Zb$ for untwisted (Neveu-Schwarz) representations and $m,n \in \Zb+\frac12$ for twisted (Ramond) representations. 

Denote by $S(\h)$ the symmetric algebra of $\h$ in $\svect(\Cb)$, the category of finite-dimensional complex super-vector spaces. Note that as a vector space, $S(\h)$ is simply the exterior algebra of the vector space underlying $\h$; in particular, $S(\h)$ is finite-dimensional. The categories of untwisted and twisted representations of $\hat\h$ (of a certain type) are equivalent to \cite[Thms.\,2.4\,\&\,2.8]{Runkel:2012cf}:
\begin{equation}
\text{(untwisted)} ~~ \Cc_0 := \Rep_{\svect} S(\h) 
\hspace{3em}
\text{(twisted)} ~~ \Cc_1 := \svect(\Cb)  \ .
\end{equation}
We would like to stress that for $\dim \h = 2n > 0$, $\Cc_0$ is not semi-simple.

A conformal field theory calculation endows the category $\Cc = \Cc_0 + \Cc_1$ with a $\Zb/2\Zb$-graded tensor product \cite[Thm.\,3.13]{Runkel:2012cf}. This tensor product is of the form stated in Section \ref{intro} with symmetric category $\Sc = \svect(\Cb)$ and Hopf algebra $H = S(\h)$. 

The associativity isomorphism is determined by the copairing $C \in \h \otimes \h$ dual to the super-symmetric pairing $(-,-)$ on $\h$, and by a top-form $\lambda$ on $S(\h)$ such that $(\lambda \otimes \lambda)(e^{-C}) = 1$ \cite[Thm.\,6.2]{Runkel:2012cf}. To be more specific, pick a basis $\{ e_i \}_{i=1,\dots,2n}$ of $\h$ such that the pairing takes the standard form $(e_{2k-1},e_{2k})=1 = -(e_{2k},e_{2k-1})$ for $k=1,\dots,n$. Then $C = \sum_{k=1}^{n} (e_{2k} \otimes e_{2k-1} - e_{2k-1} \otimes e_{2k})$ and, if we set $\hat C = -2 \sum_{k=1}^n e_{2k} \otimes e_{2k-1}$, the top-form $\lambda$ is determined by $\lambda(\hat C^n) =  n! (-2i)^n$. The explicit form of the associativity isomorphisms will be given as a special case of Theorem \ref{hopfassoc} below.

To describe the braiding, denote by $\omega_V$ the parity involution on a super-vector space $V$, and by $s_{V,W} : V \otimes W \to W \otimes V$ the symmetric structure on $\svect(\Cb)$. Then \cite[Thm.\,6.4]{Runkel:2012cf}: 
\begin{equation}\label{sf-braid}
\raisebox{.7em}{
\begin{tabular}{ccl}
 $A$ & $B$ & $c_{A,B} ~:~ A \ast B \to B \ast A$
\\
\hline
 $\Cc_0$ & $\Cc_0$ & $s_{A,B} \circ \exp(-C)$
 \\
 $\Cc_0$ & $\Cc_1$ & $s_{A,B} \circ \big( \exp(\frac12 \hat C) \otimes \id_B \big)$
\\
 $\Cc_1$ & $\Cc_0$ & $s_{A,B} \circ \big( \id_A \otimes\exp(\frac12 \hat C) \big) \circ ( \id_A \otimes \omega_B )$
\\
 $\Cc_1$ & $\Cc_1$ & $e^{-i \pi  \frac{n}4} \cdot
  \big(\id_{S(\h)} \otimes s_{A,B}\big) \circ \big(\exp( - \frac12 \hat C) \otimes \id_A \otimes  \omega_B \big)$ 
\end{tabular}}
\end{equation}
In the last line, note that for $A,B \in \Cc_1$ we have $A \ast B = S(\h) \otimes A \otimes B$.

\section{A unified framework}\label{unif}

The braiding isomorphisms \eqref{ty-braid} and \eqref{sf-braid} in the two examples just discussed may look quite different at first glance, but -- just as was the case for the tensor product $\ast$ itself -- they are actually two instances of the same structure.

\begin{figure}[bt]
\begin{center}
\begin{tabular}{c@{\hspace{0.9em}}|@{\hspace{0.9em}}c@{\hspace{0.9em}}|@{\hspace{0.9em}}c@{\hspace{0.9em}}|@{\hspace{0.9em}}c@{\hspace{0.9em}}|@{\hspace{0.9em}}c}
\small 010 & \small 011 & \small 101 & \small 110 & \small 111 \\
\hline
&&&& \\[-.5em]
\ipic{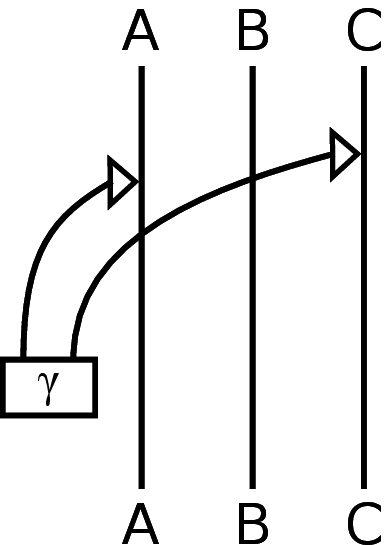} &
\ipic{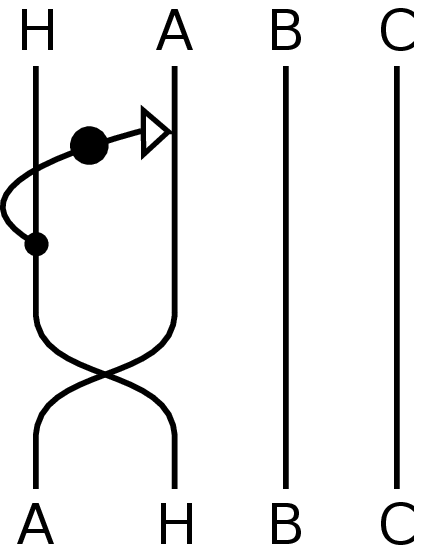} &
\ipic{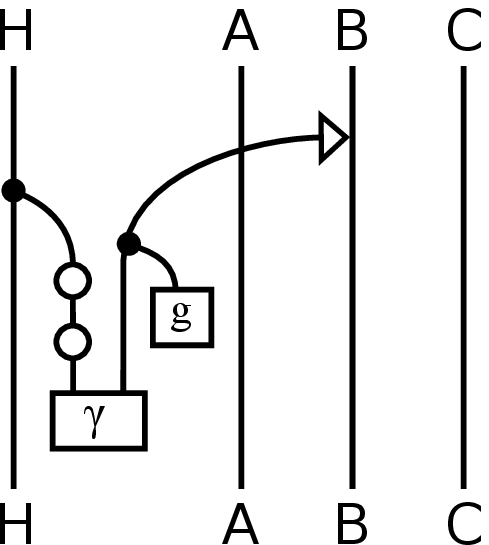} &
\ipic{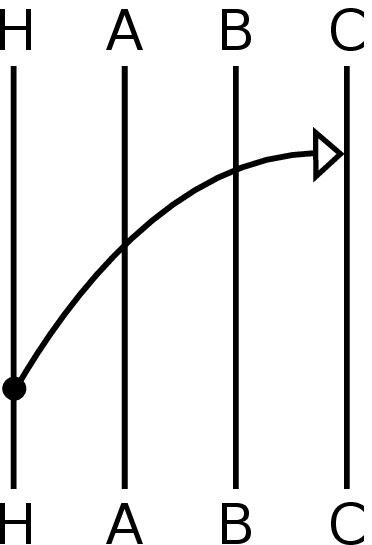} &
\ipic{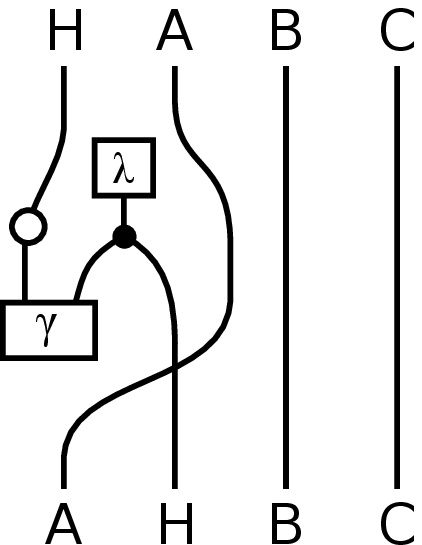} 
\end{tabular}
\end{center}
\caption{Associativity isomorphism $\alpha_{A,B,C} : A \ast (B \ast C) \to (A \ast B) \ast C$. The label $abc$ means that $A \in \Cc_a$, $B \in \Cc_b$, $C \in \Cc_c$. In the three non-listed cases 000, 001, 100, $\alpha_{A,B,C}$ is the identity (or rather the associator of the underlying category $\Sc$). The diagrams are read from bottom to top, the empty and solid dot denote $S$ and $S^{-1}$, respectively, and the three-valent vertices are the product and coproduct. The arrowhead depicts the action of $H$ on a module.}
\label{assoc}
\end{figure}

Namely, let $\Sc$ be a pivotal symmetric monoidal category (i.e.\ a ribbon category with symmetric braiding $s_{A,B} : A \otimes B \to B \otimes A$) and let $H$ be a Hopf algebra in $\Sc$ with invertible antipode. Denote by $\Rep_{\Sc}(H)$ the monoidal category of left $H$-modules in $\Sc$ and set
\begin{equation}
 \Cc = \Cc_0 + \Cc_1
 \qquad \text{with} \qquad
 \Cc_0 = \Rep_{\Sc}(H) \quad , \quad
 \Cc_1 = \Sc \ .
\end{equation} 
On $\Cc$ we fix the $\Zb/2\Zb$-graded tensor product functor $\ast$ from Section \ref{intro}. Denote by $\mu$ the multiplication of $H$, by $\Delta$ the coproduct, and by  $S$ the antipode.
For a morphisms $x : \one \to H$ in $\Sc$, write ${}_xM$ and $M_x$ for the endomorphism of $H$ given by left- and right-multiplication with $x$, respectively, and  $\Ad_x$ for the endomorphism given by $x(-)x^{-1}$.
Associativity isomorphisms for $\ast$ can be obtained from Hopf-algebraic data as follows \cite[Cor.\,3.17]{Davydov:2012xg}:

\begin{theorem}\label{hopfassoc}
Let $\gamma : \one \to H \otimes H$ and $\lambda : H \to \one$ be two morphisms in $\Sc$ such that
\begin{enumerate}
\item
$\gamma$ is a non-degenerate Hopf-copairing,
\item
$\lambda$ is a right cointegral for $H$, such that there exits $g : \one \to H$ with $(\id \otimes \lambda) \circ \Delta = g \circ \lambda$, and such that $(\lambda \otimes \lambda) \circ (\id \otimes S) \circ \gamma = \id_\one$,
\item
$\gamma$ 
satisfies the symmetry condition
$s_{H,H} \circ \gamma = (\id \otimes (S^2 \circ \Ad_g^{-1})) \circ \gamma$.
\end{enumerate}
Then the natural isomorphisms in Figure \ref{assoc} 
define associativity isomorphisms for $\ast$.
\end{theorem}

Note that it is not claimed that the above description gives {\em all} associativity isomorphisms for $\ast$; outside of $\Sc=\vect(k)$ this would require extra assumptions. 
In \cite{Davydov:2012xg}, the above theorem is actually proved in the more general setting of $\Sc$ being ribbon but not necessarily symmetric. The relevant Hopf algebra notions are reviewed in \cite[Sec.\,2]{Davydov:2012xg}. 

\begin{example}
1.\ The Tambara-Yamagami categories are recovered for $\Sc = \vect(\Cb)$ and $H = \mathrm{Fun}(G,\Cb)$. The cointegral and copairing are $\lambda = \tau \sum_{a \in G} \delta_a$ and $\gamma = \sum_{a,b \in G} \chi(a,b) \, \delta_a \otimes \delta_b$. 
\\[.3em]
2.\ For symplectic fermions take $\Sc = \svect(\Cb)$ and $H = S(\h)$. The copairing is $\gamma = e^{C}$ and the cointegral $\lambda$ is as given in Section \ref{sf}.
\\[.3em]
3.\ Another example for $\Sc = \vect(\Cb)$ is provided by Sweedler's four-dimensional Hopf algebra, which is not semi-simple \cite[Sec.\,3.8.3]{Davydov:2012xg}. This illustrates that Theorem \ref{hopfassoc} is more general than Tambara-Yamagami categories even in the vector space case.
\end{example}

\begin{figure}[bt]
\begin{center}
\begin{tabular}{c@{\hspace{1.4em}}|@{\hspace{1.4em}}c@{\hspace{1.4em}}|@{\hspace{1.4em}}c@{\hspace{1.4em}}|@{\hspace{1.4em}}c}
\small 00 & \small 01 & \small 10 & \small 11  \\
\hline
&&& \\[-.5em]
\jpic{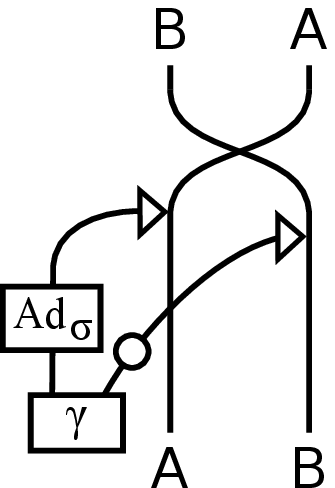} &
\jpic{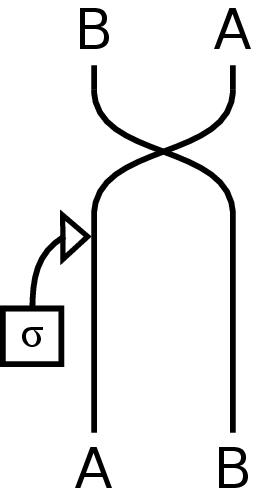} &
\jpic{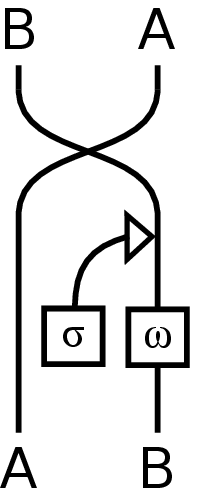} &
$\beta ~\cdot~$ \jpic{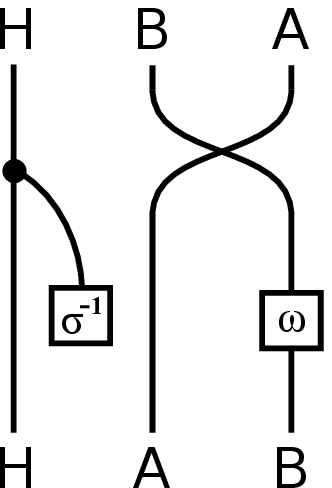} 
\end{tabular}
\end{center}
\caption{Braiding isomorphisms $c_{A,B} : A \ast B \to B \ast A$. The notation is as in Figure \ref{assoc}.}
\label{braiding}
\end{figure}

For the braiding isomorphisms we need to fix an involutive monoidal automorphism $\omega$ of the identity functor on $\Sc$. For $\Sc = \vect(\Cb)$, $\omega$ is necessarily the identity, but for $\Sc = \svect(\Cb)$ there are already two choices: the identity and parity involution. We have \cite[Thm.\,1.2\,\&\,Rem.\,4.11]{Davydov:2012xg}:

\begin{theorem} \label{hopfbraid}
Let $H$, $\gamma$, $\lambda$ and $g$ be as in Theorem \ref{hopfassoc}
and let $\sigma : \one \to H$ and $\beta : \one \to \one$ be invertible.
Suppose that
\begin{enumerate}
\item 
$\gamma$ is determined through $\sigma$ by 
$\gamma =  ({}_{\sigma^{-1}}M \otimes M_{\sigma^{-1}}) \circ \Delta \circ \sigma$.
\item 
$\lambda$ satisfies $\lambda \circ S = \lambda \circ \Ad_\sigma$
and $\lambda \circ \sigma = \beta \circ \beta$.
\item
$\Ad_\sigma$ is a Hopf-algebra isomorphism $H \to H_\mathrm{cop}$ (the opposite coalgebra). 
\item
$S \circ \sigma = {}_gM \circ \sigma = M_{g^{-1}} \circ \sigma$.
\item
$\omega$ evaluated on $H$ satisfies
$(\id \otimes \omega_H) \circ \gamma = \big(  \Ad_\sigma \otimes (\Ad_\sigma^{-1}\circ S) \big) \circ \gamma$.
\end{enumerate}
Then the natural isomorphisms in Figure \ref{braiding} define a braiding on $\Cc$.
\end{theorem}

If $\sigma^2$ is central in $H$, then $\Cc$ can be made into a ribbon category with twist isomorphisms $\theta_A = \sigma^{-2}.(-)$ for $A \in \Cc_0$ (i.e.\ the left action of $\sigma^{-2}$ on the $H$-module $A$), and $\theta_A = \beta^{-1} \, \omega_A$ for $A \in \Cc_1$, see \cite[Prop.\,4.18]{Davydov:2012xg}.

\begin{example} 
1.\ In the Tambara-Yamagami case, $\sigma$ and $\beta$ are as in Section \ref{ty}. Comparing \eqref{ty-braid} and Figure \ref{braiding} shows that the braiding isomorphisms match ($S=\id_H$ for an elementary 2-group, and $\Ad_\sigma=\id_H$ since $H$ is commutative). 
\\[.3em]
2.\ For symplectic fermions choose $\sigma = \exp(\frac12 \hat C)$ and $\beta = e^{- \pi i n / 4}$. Then \eqref{sf-braid} agrees with Figure \ref{braiding}.
\\[.3em]
3.\ Sweedler's Hopf algebra is quasi-triangular, but the resulting braiding on $\Cc_0$ does not extend to all of $\Cc$ (at least via the above construction). However, one can find a 16-dimensional semi-simple Hopf algebra in $\vect(\Cb)$ which is neither commutative nor co-commutative for which Theorems \ref{hopfassoc} and \ref{hopfbraid} apply \cite[Sec.\,4.7.4]{Davydov:2012xg}. This is another instance where our setting is more general than the Tambara-Yamagami case even for $\Sc=\vect(\Cb)$.
\end{example}

\newcommand\arxiv[2]      {\href{http://arXiv.org/abs/#1}{#2}}
\newcommand\doi[2]        {\href{http://dx.doi.org/#1}{#2}}
\newcommand\httpurl[2]    {\href{http://#1}{#2}}

\end{document}